\documentclass[%
aip,
cp,  % Conference Proceedings
amsmath,amssymb,%nobibnotes,
reprint,%
]{revtex4-2}

\usepackage{graphicx}% Include figure files
\usepackage{dcolumn}% Align table columns on decimal point
\usepackage{bm}% bold math
\usepackage[mathlines]{lineno}% Enable numbering of text and display math
\usepackage{subcaption}
\usepackage[utf8]{inputenc}
\usepackage[T1]{fontenc}
\usepackage{mathptmx} 
\usepackage{amsthm}
\newtheorem{theorem}{Theorem} % Lingkungan teorema
\begin{document}
	
	\title{Dynamical Analysis of a Predator-Prey Model with Additif Allee Effect and Prey Group Defense}% Force line breaks with \\
	
	\author{Resmawan} % Write as First name Surname
	\email[Corresponding author: ]{resmawan@student.ub.ac.id}
	\author{Agus Suryanto}%
	\email{suryanto@ub.ac.id}
	\author{Isnani Darti}%
	\email{isnanidarti@ub.ac.id}
	\affiliation{
		Department of Mathematics, Faculty of Mathematics and Natural Sciences, University of Brawijaya, Indonesia
	}
	
	\author{Hasan S Panigoro}
	\email{hspanigoro@ung.ac.id}
	\affiliation{%
		Department of Mathematics, Faculty of Mathematics and Natural Sciences, Universitas Negeri Gorontalo, Indonesia
	}
	
	The International Symposium on Biomathematics (Symomath), 10-12 July 2024, Depok, Indonesia
	
	\begin{abstract}
		In this article, we develop a predator-prey model with Allee effect and prey group defense. The model has three equilibrium points i.e. the trivial point, the predator extinction point, and the coexistence point. All equilibrium points are locally asymptotically stable under certain conditions. The Allee effect in this model influences the stability of the equilibrium point. If the Allee effect is weak, then the trivial equilibrium point is unstable. Meanwhile, if the Allee effect is strong, then the trivial equilibrium point is locally asymptotically stable. Those mean that a strong Allee effect can lead to the extinction of both populations. Moreover, under weak Allee condition, forward bifurcation and Hopf bifurcation occur at the predator extinction equilibrium point. Meanwhile, a strong Allee effect may induce bistability at both the trivial equilibrium point and the predator extinction equilibrium point. Those mean that prey can survive without the presence of predators, but a strong Allee effect can lead to prey extinction if the population size is very small. To support our analytical findings, we perform some numerical simulations in the final section.
	\end{abstract}
	
	\maketitle
	
	\section{\label{sec:level1}Introduction}
	
	Mathematical models representing predator-prey interactions are extensively studied in applied mathematics. Lotka \cite{Lotka1925} and Volterra \cite{Volterra1927} introduced models of interactions between two species, adopting the population growth model from \cite{Malthus1872}. This model, later known as the Lotka-Volterra model, has become a reference for the development of other foundational models, such as those by \cite{Leslie1960} and \cite{Rosenzweig1963}. The discussion of predator-prey models remains crucial for studying various biological phenomena. Efforts to construct these models continue in order to develop models that are more realistic and consistent with observed biological phenomena.
	
	The development of models has extensively considered various biological phenomena observed in both prey and predators. Some models have been developed by taking into account specific phenomena such as differences in the age structure of prey \cite{Beay2020, Zhang2022}, the effects of fear on prey \cite{Purnomo2023}, internal competition among the same species \cite{Long2022, Panigoro2023}, and the presence of disease in species \cite{Purnomo2023}. Predator-prey models are always open to further examination with additional considerations in line with ecological factors observed in specific species.
	
	A crucial ecological factor, especially for species facing extinction, is the Allee effect, which describes the difficulties in regeneration leading to extinction threats. Predator-prey models incorporating the Allee effect have been widely discussed, examining its influence on classical models such as the Leslie-Gower model \cite{Rahmi2021a}, Lotka-Volterra model \cite{Panigoro2021}, and Rosenzweig-MacArthur model \cite{Panigoro2022b}. Studies also explore the Allee effect with different functional responses such as Beddington-DeAngelis \cite{Rahmi2021b} and Michaelis-Menten \cite{Panigoro2022, Rahmi2022}. Additionally, the Allee effect combined with intraspecific competition in predators is discussed by \cite{Anggriani2023}, and its impact on eco-epidemiological models by \cite{Sidik2022, Rahmi2023}. Models considering the Allee effect on predators can be found in the studies by \cite{Rahmi2021a, Panigoro2022, Anggriani2023}, while the consideration of the Allee effect on prey is discussed in the models by \cite{Ye2021, Bai2022, Xie2023, GarciaCuenca2023, Xie2023}.
	
	Besides the Allee effect, prey group defense against predator attacks is another observed phenomenon \cite{Zhang2021}. Several models have considered this phenomenon. Zhang et al. \cite{Zhang2021} discuss Hopf bifurcation in predator-prey models with prey group defense and time delays, representing group defense as an exponential function. Jiao et al. \cite{JIAO2022} develop a Leslie-Gower model incorporating prey group defense with a threshold value, using a type IV functional response. Lynch \cite{Lynch2017} also considers prey group defense using a modified Holling type IV functional response. Prey group defense behavior can influence predator population density and enhance prey survival.
	
	This paper is organized as follows. The model structure is described in Section 2, followed by the existence and local stability analysis of the equilibrium points of the model are discussed in Section 3, while numerical simulations and interpretations are presented in Section 4. Finally, we draw some concluding remarks in Section 5.
	
	\section{\label{sec:level1}MODEL}
	Developing models that consider biological phenomena in both predator and prey species is an ongoing area for improvement to create more realistic models. In this study, we develop a predator-prey model that incorporates the Allee effect and prey group defense. We reference the model by \cite{Bai2022} to include the Allee effect and the model by \cite{JIAO2022} to account for prey group defense. The predator-prey model by \cite{Bai2022} combines the Allee effect with a Holling type I functional response. We modify this model by using a Holling type IV functional response, which is more ecologically relevant and represents the prey group defense. This response shows that larger prey groups experience lower predation rates. While \cite{JIAO2022} use the Holling type IV functional response to represent prey group defense, they do not consider the Allee effect. In our study, we construct a model that incorporates the Holling type IV functional response for prey group defense from \cite{JIAO2022}, along with the Allee effect on the prey as discussed by \cite{Bai2022}. The constructed model involves two variables: prey population density $(N)$ and predator population density $(P)$ at time $t$. This model is referred to as the predator-prey model with the Allee effect and prey group defense, presented in equation \eqref{GrindEQ__2_},
	\begin{eqnarray}\label{GrindEQ__2_}
		\frac{dN}{dt}=&&\ rN\left(1-\frac{N}{K}-\frac{h}{w+N}\right)-\frac{aNP}{b+N^2},\\ 
		\frac{dP}{dt}=&&\ \frac{cNP}{b+N^2}-\delta P,\nonumber
	\end{eqnarray}
	with $r$ is the intrinsic growth rates of the prey. Furthermore, the positive constants $a,b, \delta$ and $K$ represent the prey's predation rate, the half-saturation constant of predation, the death rate of the predator, and the prey's carrying capacity. Meanwhile, $h, w > 0$ describes the degree of the Allee effect, with $h$ the rate for severity of Allee, and $w < K$ the prey population size at which fitness is half of its maximum value. In particular, if $h < w $ or $w < h$, then the system's \eqref{GrindEQ__2_} has weak or strong Allee effect, respectively \cite{GarciaCuenca2023}. 
	
	\section{Existence and Local Stability of Equilibrium Points}
	\subsection{The Existence of Equilibrium Points }
	The equilibrium points of model \eqref{GrindEQ__2_} are obtained by simultaneously solving $\frac{dN}{dt}=0$ and $\frac{dP}{dt}=0,\ $i.e.
	\begin{eqnarray}\label{GrindEQ__5_}
		rN\left(1-\frac{N}{K}-\frac{h}{w+N}\right)-\frac{aNP}{b+N^2}=&&\ 0 \\ 
		\frac{cNP}{b+N^2}-\delta P=&&\ 0.\nonumber
	\end{eqnarray} 
	By solving system \eqref{GrindEQ__5_}, three types of equilibrium points are obtained: trivial equilibrium points, axial equilibrium points, and coexistence equilibrium points. The trivial equilibrium point represents the condition of extinction for all populations. The trivial equilibrium point is denoted by $E_0=\left(0,0\right)$, which always exists in $R^2_+$. 
	
	The axial equilibrium points are denoted by $E_n=\left(N_n,0\right),\ n=1,2,3$, representing the equilibrium points where the predator population is extinct, with $N$ obtained from equation \eqref{GrindEQ__6_}:
	\begin{equation}\label{GrindEQ__6_}
		N^2-\left(K-w\right)N+K\left(h-w\right)=0.
	\end{equation} 
	Suppose $N_1$ and $N_2$ are the two roots of equation \eqref{GrindEQ__6_}, then it follows:
	\begin{eqnarray}
		N_1=&&\ \frac{\left(K-w\right)+\sqrt{{\left(K-w\right)}^2-4K\left(h-w\right)}}{2}\label{GrindEQ__7_} \\
		N_2=&&\ \frac{\left(K-w\right)-\sqrt{{\left(K-w\right)}^2-4K\left(h-w\right)}}{2}.\label{GrindEQ__8_}
	\end{eqnarray} 
	
	The existence of $N_n$ can be determined by examining the Allee effect condition $(h-w)$, the value of $(K-w)$, and the discriminant value from equation \eqref{GrindEQ__6_}, namely:
	\begin{equation}\label{GrindEQ__9_}
		D_1={\left(K-w\right)}^2-4K\left(h-w\right).
	\end{equation}
	
	\begin{itemize}
		\item  Weak Allee Effect Case\\
		The weak Allee effect in system \eqref{GrindEQ__2_} occurs if $h\ <\ w$. If $h\ <\ w$, then $D_1\ >\ 0$, which leads to the existence of axial equilibrium points depending on the value of $K\ -\ w$:
		\begin{enumerate}
			\item [a.] If $K>w,\ \ $then $\sqrt{{\left(K-w\right)}^2-4K\left(h-w\right)}>(K-w)$, so that $N_1>0$ and $N_2<0.$ 
			
			\item [b.] If $K<w,$ then $\sqrt{{\left(K-w\right)}^2-4K\left(h-w\right)}=\mathrm{\ }\sqrt{{\left(w-K\right)}^2-4K\left(h-w\right)}>\left(w-K\right)$, so that $N_1>0$ and $N_2<0.$
			
			Hence, if there is a weak Allee effect in the system \eqref{GrindEQ__2_}, then there exists one axial equilibrium point, namely $E_1=\left(N_1,0\right)$.
		\end{enumerate}
		
		\item  Strong Allee Effect Case \\
		The strong Allee effect in system \eqref{GrindEQ__2_} occurs if $h>\ w$. If $h>\ w$, then the existence of axial equilibrium points depends on the discriminant $\left(D_1\right)$ in Eq \eqref{GrindEQ__9_} and $(K-w)$ value:
		
		\begin{enumerate}
			\item [a.] $D<0$ case,\\
			If $D<0,$ then the axial equilibrium point $E_n=\left(N_n,0\right)$ do not exist.
			
			\item [b.] $D>0$ case, \\
			If $D>0,$ then $\ h<\frac{{\left(K+w\right)}^2}{4K}.$ Furthermore, the existence of the axial equilibrium point depends on $\left(K-w\right)$ value:
			
			\begin{enumerate}
				\item [(i)]  If $K>w,$ then $\sqrt{{\left(K-w\right)}^2-4K\left(h-w\right)}<(K-w)$, so that $N_1>0$ and $N_2>0.$ In this case, there are two existing axial equilibrium points, namely $E_1=\left(N_1,0\right)$ and $E_2=\left(N_2,0\right).$
				
				\item [(ii)] If $K<w,$ then $\sqrt{{\left(K-w\right)}^2-4K\left(h-w\right)}<\left(K-w\right),$ so that $N_1<0$ and $N_2<0.$
				In this case, $E_1=\left(N_1,0\right)$ and $E_2=\left(N_2,0\right)$ do not exist.
			\end{enumerate}
			
			\item [c.] $D=0$ case,\\
			If $D=0,$ then there exists one existing axial equilibrium point, i.e $E_3=\left(N_3,0\right),$ with
			\[N_3=\frac{K-w}{2}.\] 
			$E_3$ exists if $K>w$ and do not exist if $K<w$.
		\end{enumerate}
	\end{itemize}
	
	Therefore, there are three axial equilibrium points, namely $E_1=\left(N_1,0\right)$, $E_2=\left(N_2,0\right)$, and $E_3=\left(N_3,0\right)$, whose existency depends on the Allee effect conditions. The existence of axial equilibrium points is stated in Theorem \ref{teo4} and Theorem \ref{teo5}.
	
	\begin{theorem}\label{teo4}
		If the system \eqref{GrindEQ__2_} has a weak Allee effect $(h<w)$, then the equilibrium point $E_1=\left(N_1,0\right)$ exists and is unique.
	\end{theorem}
	
	\begin{theorem}\label{teo5}
		Let $K>w$ and the system \eqref{GrindEQ__2_} has a strong Allee effect $(h>w)$:
		
		\begin{enumerate}
			\item [(i)] If $h>\frac{{\left(K+w\right)}^2}{4K},$ then there are no axial equilibrium points.
			
			\item [(ii)] If $h=\frac{{\left(K+w\right)}^2}{4K},$ then there exists exactly one axial equilibrium point, namely $E_3$.
			
			\item [(iii)] If $h<\frac{{\left(K+w\right)}^2}{4K},$ then there are two axial equilibrium point, namely $E_1$ and $E_2.$
		\end{enumerate}
	\end{theorem}
	
	The coexistence equilibrium points are denoted by $E_i=\left(N_i,P_i\right),\ i=4,5,6$, which represent the condition where all populations exist, with $N_i$ and $P_i$ obtained from Eq. \eqref{GrindEQ__10_} and \eqref{GrindEQ__11_}.
	\begin{eqnarray}
		\frac{cN_i}{b+N^2_i}-\delta =&&\ 0\label{GrindEQ__10_}\\ r\left(1-\frac{N_i}{K}-\frac{h}{w+N_i}\right)-\frac{aP_i}{b+N^2_i}=&&\ 0.\label{GrindEQ__11_}
	\end{eqnarray} 
	
	\begin{itemize}
		\item  From Eq. \eqref{GrindEQ__10_}, we obtain
		\begin{equation}\label{GrindEQ__12_}
			D_2=c^2-4b{\delta }^2,\ \ N_4=\frac{c+\sqrt{D_2}}{2\delta }\ \ \mathrm{and}\ \ N_5=\frac{c-\sqrt{D_2}}{2\delta }.
		\end{equation} 
		$N_{4,5}$ exists if $D_2\ge 0$ or $b\le {\left(\frac{c}{2\delta }\right)}^2.$
		
		\item From Eq. \eqref{GrindEQ__11_}, we obtain \\ 
		\begin{equation}\label{GrindEQ__13_}
			P_i(N_i)=\frac{r\left(b+N^2_i\right)\left[\left(K-w\right)N_i-K\left(h-w\right)-N^2_i\right]}{Ka(w+N_i)}, i=4,5,6
		\end{equation}
	\end{itemize}
	
	The existence of the coexistence equilibrium point is stated in Theorem \ref{teo6}.
	
	\begin{theorem}\label{teo6}
		Let $D_2=c^2-4b{\delta }^2,$ $N_4=\frac{c+\sqrt{D_2}}{2\delta },N_5=\frac{c-\sqrt{D_2}}{2\delta },\ N_6=\frac{c}{2\delta },$ and $P_i\left(N_i\right)=\frac{r\left(b+N^2_i\right)\left[\left(K-w\right)N_i-K\left(h-w\right)-N^2_i\right]}{Ka(w+N_i)}, i=4,5,6$:
		\begin{enumerate}
			\item [(i)] If  $b>{\left(\frac{c}{2\delta }\right)}^2,$ then there are no coexistence equilibrium points.
			
			\item [(ii)] If  $b={\left(\frac{c}{2\delta }\right)}^2$and $P\left(N_6\right)>0,$ then there exists excacly one coexistence equilibrium point, namely $E_6=\left(N_6,P\left(N_6\right)\right).$
			
			\item [(iii)] If $b<{\left(\frac{c}{2\delta }\right)}^2,P\left(N_4\right)>0,$ and $P\left(N_5\right)>0,$ then there are two coexistence equilibrium points, namely $E_4=\left(N_4,P\left(N_4\right)\right)$ and $E_5=\left(N_5,P\left(N_5\right)\right).$
		\end{enumerate}
	\end{theorem} 
	
	\subsection{Local Stability}
	Linearization around the equilibrium point is carried out so that the Jacobian matrix is obtained:
	\begin{equation}\label{GrindEQ__14_}
		J=\left[ \begin{array}{cc}
			r-\frac{2r}{K}N-\frac{rhw}{{\left(w+N\right)}^2}-\frac{aP\left(b-N^2\right)}{{\left(b+N^2\right)}^2} & -\frac{aN}{b+N^2} \\ 
			\frac{cP}{b+N^2}-\frac{2cN^2P}{{\left(b+N^2\right)}^2} & \frac{cN}{b+N^2}-\delta  \end{array}
		\right].
	\end{equation}

	By substituting $E_0=\left(0,0\right)$ to Eq. \eqref{GrindEQ__14_}, we obtain
	\begin{equation*}
		J_{E_0}=\left[ \begin{array}{cc}
			\frac{r\left(w-h\right)}{w} & 0 \\ 
			0 & -\delta  \end{array}
		\right]
	\end{equation*}
	and we get two eigen values ${\lambda }_1=\frac{r\left(w-h\right)}{w}$ and ${\lambda }_2=-\delta <0.$ Hence, $E_0$ is locally asymptotically stable if $h>w$ and unstable (saddle node) if $h<w$. The stability of trivial equilibrium point is stated in Theorem \ref{teo7}.
	
	\begin{theorem}\label{teo7}
		The trivial equilibrium point, $E_0=\left(0,0\right)$, is locally asymptotically stable if the Allee effect is strong $\left(h>w\right)\ $and unstable (saddle-node) if the Allee effect is weak $\left(h<w\right).$
	\end{theorem}
	
	%\subsubsection{Stability Analysis of Axial Equilibrium Points}
	By substituting $E_1=\left(N_1,0\right)$ to Eq. \eqref{GrindEQ__14_}, we obtain 
	\begin{equation*}
		J_{E_1}=\left[ \begin{array}{cc}
			rN_1\left(\frac{h}{{\left(w+N_1\right)}^2}-\frac{1}{K}\right) & -\frac{aN_1}{b+N^2_1} \\ 
			0 & \frac{cN_1-\delta \left(b+N^2_1\right)}{b+{N_1}^2} \end{array}
		\right].
	\end{equation*} 
	The eigen values of $J_{E_1}$ are ${\lambda }_1=rN_1\left(\frac{h}{{\left(w+N_1\right)}^2}-\frac{1}{K}\right)\ \mathrm{and}\mathrm{\ }{\lambda }_2=\frac{cN_1-\delta \left(b+N^2_1\right)}{b+{N_1}^2}.$ If the Allee effect is weak, then $h<\frac{{\left(K+w\right)}^2}{4K},$ and it can be shown that ${\lambda }_1<0.$ Furthermore, it can be shown that the value of ${\lambda }_2$ depends on $c$. If $c<\frac{\delta \left(b+N^2_1\right)}{N_1},\ $then ${\lambda }_2<0,$ making the axial equilibrium point, $E_1,$ is locally asymptotically stable and if $c>\frac{\delta \left(b+N^2_1\right)}{N_1},\ $then ${\lambda }_2>0,$ making the axial equilibrium point is unstable (saddle-node). The stability of axial equilibrium point with a weak Allee effect is stated in Theorem \ref{teo8}.
	
	\begin{theorem}\label{teo8}
		Suppose the system \eqref{GrindEQ__2_} has a weak Allee effect. The axial equilibrium point, $E_1=\left(N_1,0\right)$ is locally asymptotically stable if $c<\frac{\delta \left(b+N^2_1\right)}{N_1}$ and unstable (saddle node) if  $c>\frac{\delta \left(b+N^2_1\right)}{N_1}.$
	\end{theorem}
	
	The Jacobian matrix of $E_n=(N,0)$ is 
	\begin{equation*}
		J_{E_n}=\left[ \begin{array}{cc}
			rN_n\left(\frac{h}{{\left(w+N_n\right)}^2}-\frac{1}{K}\right) & -\frac{aN_n}{b+N^2_n} \\ 
			0 & \frac{cN-\delta \left(b+N^2_n\right)}{b+N^2_n} \end{array}
		\right]
	\end{equation*} 
	where its eigenvalues are,
	\begin{equation*}
		{\lambda }_1=rN_n\left(\frac{h}{{\left(w+N_n\right)}^2}-\frac{1}{K}\right)\ \ \ \ \ \ \mathrm{and}\mathrm{\ \ \ \ \ \ }{\lambda }_2=\frac{cN_n-\delta \left(b+N^2_n\right)}{b+N^2_n}.
	\end{equation*}
	
	\noindent We can show that ${\lambda }_2<0$ if $c<\frac{\delta \left(b+N^2_n\right)}{N_n}.$ Furthermore, ${\lambda }_1$ is depend on the Allee effect case. For the strong Allee effect $(h>w)$, we have the following case:
	\begin{itemize}
		\item  If $h=\frac{{\left(K+w\right)}^2}{4K},$ then
		\begin{equation*}
			{\lambda }_1=rN_3\left(\frac{h}{{\left(w+N_3\right)}^2}-\frac{1}{K}\right)=rN_3\left(\frac{\frac{{\left(K+w\right)}^2}{4K}}{{\left(\frac{K+w}{2}\right)}^2}-\frac{1}{K}\right)=0.
		\end{equation*}
		Since ${\lambda }_1=0$, the axial equilibrium point $E_3=\left(\frac{K-w}{2},0\right)$ is non-hyperbolic. 
		
		\item  If  $h<\frac{{\left(K+w\right)}^2}{4K},$ then
		\begin{eqnarray*}
			{\lambda }_1=&&\ rN_n\left(\frac{h}{{\left(w+N_n\right)}^2}-\frac{1}{K}\right)<\ \frac{rN_{1,2}}{K}\left(\frac{{\left(K+w\right)}^2-{\left(\left(K+w\right)\pm \sqrt{{\left(K+w\right)}^2-4Kh}\right)}^2}{{\left(\left(K+w\right)\pm \sqrt{{\left(K+w\right)}^2-4Kh}\right)}^2}\right).
		\end{eqnarray*}
		
		\begin{enumerate}
			\item [(i)] For $E_1,$
			\begin{equation*}
				{\lambda }_1<\frac{rN_1}{K}\left(\frac{{\left(K+w\right)}^2-{\left(\left(K+w\right)+\sqrt{{\left(K+w\right)}^2-4Kh}\right)}^2}{{\left(\left(K+w\right)+\sqrt{{\left(K+w\right)}^2-4Kh}\right)}^2}\right)
			\end{equation*} 
			Furthermore, it is shown that 
			\begin{eqnarray*}
				{\left(K+w\right)}^2-{\left(\left(K+w\right)+\sqrt{{\left(K+w\right)}^2-4Kh}\right)}^2=&&-2\left(K+w\right)\sqrt{{\left(K+w\right)}^2-4Kh}-{\left(K+w\right)}^2+4Kh\\ 
				<&&\ -2\left(K+w\right)\sqrt{{\left(K+w\right)}^2-4Kh}-{\left(K+w\right)}^2+{\left(K+w\right)}^2\\=&&\ -2\left(K+w\right)\sqrt{{\left(K+w\right)}^2-4Kh}<0.
			\end{eqnarray*} 
			Since ${\lambda }_1<0$ and ${\lambda }_2<0,$ then $E_1$ is locally asymptotically stable.
			
			\item [(ii)] For $E_2,$
			\begin{equation*}
				{\lambda }_1<\frac{rN_2}{K}\left(\frac{{\left(K+w\right)}^2-{\left(\left(K+w\right)-\sqrt{{\left(K+w\right)}^2-4Kh}\right)}^2}{{\left(\left(K+w\right)-\sqrt{{\left(K+w\right)}^2-4Kh}\right)}^2}\right)
			\end{equation*} 
			Furthermore, it is shown that
			\begin{eqnarray*}
				{\left(K+w\right)}^2-{\left(\left(K+w\right)-\sqrt{{\left(K+w\right)}^2-4Kh}\right)}^2=&&\ 2\left(K+w\right)\sqrt{{\left(K+w\right)}^2-4Kh}-{\left(K+w\right)}^2+4Kh\\ 
				<&&\ 2\left(K+w\right)\sqrt{{\left(K+w\right)}^2-4Kh}-{\left(K+w\right)}^2+{\left(K+w\right)}^2\\ 
				=&&\ 2\left(K+w\right)\sqrt{{\left(K+w\right)}^2-4Kh}>0.
			\end{eqnarray*} 
			Since ${\lambda }_1>0$ and ${\lambda }_2<0,$ then $E_2$ is unstable (saddle-node).
		\end{enumerate}
	\end{itemize}
	The stability of axial equilibrium point with a strong Allee effect is stated in Theorem \ref{teo9}.
	
	\begin{theorem}\label{teo9}
		Suppose $K>w,\ c<\frac{\delta \left(b+N^2\right)}{N}$ and the system \eqref{GrindEQ__2_} has a strong Allee effect:
		\begin{enumerate}
			\item [(i)] If $h=\frac{{\left(K+w\right)}^2}{4K},$ then the axial equilibrium point, $E_3$ is non-hyperbolic
			
			\item [(ii)] If $h<\frac{{\left(K+w\right)}^2}{4K},$ then the axial equilibrium point, $E_1$ is locally asymptotically stable and the axial equilibrium point, $E_2$ is unstable (saddle node). 
		\end{enumerate}
	\end{theorem}

	By substituting $E_i=\left(N_i,P_i\right)$ to Eq. \eqref{GrindEQ__14_}, we obtain
	\begin{equation}\label{GrindEQ__15_}
		J_{E_i}=\left[ \begin{array}{cc}
			\frac{hrN_i}{{\left(w+N_i\right)}^2\ }+\frac{2aN^2_iP_i}{{\left(b+N^2_i\right)}^2}-\frac{rN_i}{K} & -\frac{\delta }{c} \\ 
			\frac{\left(c-2\delta N_i\right)P_i}{b+N^2_i} & 0 \end{array}
		\right]
	\end{equation} 
	where
	\begin{eqnarray*}
		N_4=&&\ \frac{c+\sqrt{c^2-4{\delta }^2b}}{2\delta },\ \ \ N_5=\frac{c-\sqrt{c^2-4{\delta }^2b}}{2\delta },\ \ \ N_6=\frac{c}{2\delta }\\ 
		P_i\left(N_i\right)=&&\ \frac{r\left(b+N^2_i\right)\left[\left(K-w\right)N_i-K\left(h-w\right)-N^2_i\right]}{Ka(w+N_i)},\ \ \ i=4,5,6.
	\end{eqnarray*} 
	From Eq. \eqref{GrindEQ__15_}, the determinant and trace of the Jacobian matrix are obtained as follows:
	\begin{eqnarray*}
		\mathrm{det} \left(J_{E_i}\right)=&&\ -\left(-\frac{\delta }{c}\right)\left(\frac{\left(c-2\delta N_i\right)P}{b+N^2_i}\right)=\frac{\delta \left(c-2\delta N_i\right)P}{c\left(b+N^2_i\right)}\\ 
		\mathrm{tr}\left(J_{E_i}\right)=&&\ \frac{hrN_i}{{\left(w+N_i\right)}^2\ }+\frac{2aPN^2_i}{{\left(b+N^2_i\right)}^2}-\frac{rN_i}{K}.
	\end{eqnarray*}
	Furthermore, the stability of $E_i$ can be determined by examining the determinant and trace of $J_{E_i}$.
	
	\begin{itemize}
		\item  For $E_6=\left(N_6,P\left(N_6\right)\right),$
		\begin{equation*}
			\mathrm{det}\ \left(J_{E_6}\right)=\frac{\delta \left(c-2\delta N\right)P\left(N_6\right)}{c\left(b+N^2_6\right)}=\frac{\delta \left(c-c\right)P\left(N_6\right)}{c\left(b+N^2_6\right)}=0.
		\end{equation*} 
		Since $\mathrm{det}\ \left(J_{E_6}\right)=0,$ then $E_6$ is non-hyperbolic.
		
		\item  For $E_4=\left(N_4,P\left(N_4\right)\right),$ ${\mathrm{det} \left(J_{E_4}\right)}>0$ if $c-2\delta N_4>0$ holds.
		\begin{equation*}
			c-2\delta N_4=c-2\delta \left(\frac{c+\sqrt{c^2-4{\delta }^2b}}{2\delta }\right)=-\sqrt{c^2-4{\delta }^2b}<0.
		\end{equation*} 
		So, ${\mathrm{det} \left(J_{E_4}\right) }<0\ $for $E_4=\left(N_4,P\left(N_4\right)\right).$ Since ${\mathrm{det} \left(J_{E_4}\right) }<0,\ $then, $E_4$ is unstable. Moreover, if $tr\left(J_{E_4}\right)<0,$ then $E_4$ is saddle node.
		
		\item  For $E_5=\left(N_5,P\left(N_5\right)\right),$ ${\mathrm{det} \left(J_5\right)}>0$ if $c-2\delta N_5>0$ holds.
		\begin{equation*}
			c-2\delta N_5=c-2\delta \left(\frac{c-\sqrt{c^2-4{\delta }^2b}}{2\delta }\right)=\sqrt{c^2-4{\delta }^2b}>0.
		\end{equation*}
		So, ${\mathrm{det} \left(J_{E_5}\right) }>0$ for $E_5=\left(N_5,P\left(N_5\right)\right).$ Since ${\mathrm{det} \left(J_{E_4}\right)}>0$, $E_5$ is locally asymtotically stable if $tr\left(J_{E_5}\right)<0$.
	\end{itemize}
	The stability of coexistence equilibrium point is stated in Theorem \ref{teo10}.
	
	\begin{theorem}\label{teo10}
		Suppose $D=c^2-4b{\delta }^2,$ $N_4=\frac{c+\sqrt{D}}{2\delta },N_5=\frac{c-\sqrt{D}}{2\delta },\ N_6=\frac{c}{2\delta },$ and $P_i\left(N\right)=\frac{r\left(b+N^2_i\right)\left[\left(K-w\right)N_i-K\left(h-w\right)-N^2_i\right]}{Ka(w+N_i)}.$ Also, suppose $P\left(N_4\right)>0,$ $P\left(N_5\right)>0,\ $and $P\left(N_6\right)>0.$
		\begin{enumerate}
			\item [(i)] If  $b={\left(\frac{c}{2\delta }\right)}^2,$ then the coexistence equilibrium point $E_6$ is non-hyperbolic.
			
			\item [(ii)] If $b<{\left(\frac{c}{2\delta }\right)}^2,$ then the coesistence equilibrium point $E_4$ is unstable. Moreover, if $\mathrm{tr}\left(J_{E_4}\right)<0,$ then $E_4$ is saddle node and $E_5$ is locally asymptotically stable.
		\end{enumerate}
	\end{theorem}
	
	\section{Numerical Simulation}
	\subsection{The Impact of Predation Conversion Rate}
\begin{figure}[h!]
		%\centering
		\begin{subfigure}{0.497\textwidth}
			\centering
			\includegraphics[width=2.4in]{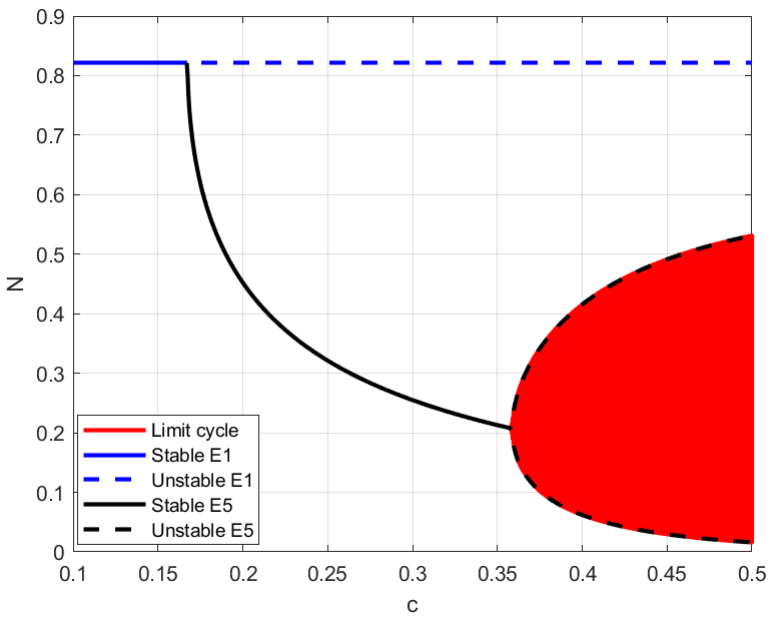}
			\caption{$ c-N $}
			\label{fig1a}
		\end{subfigure}
		\hfill
		\begin{subfigure}{0.497\textwidth}
			\centering
			\includegraphics[width=2.4in]{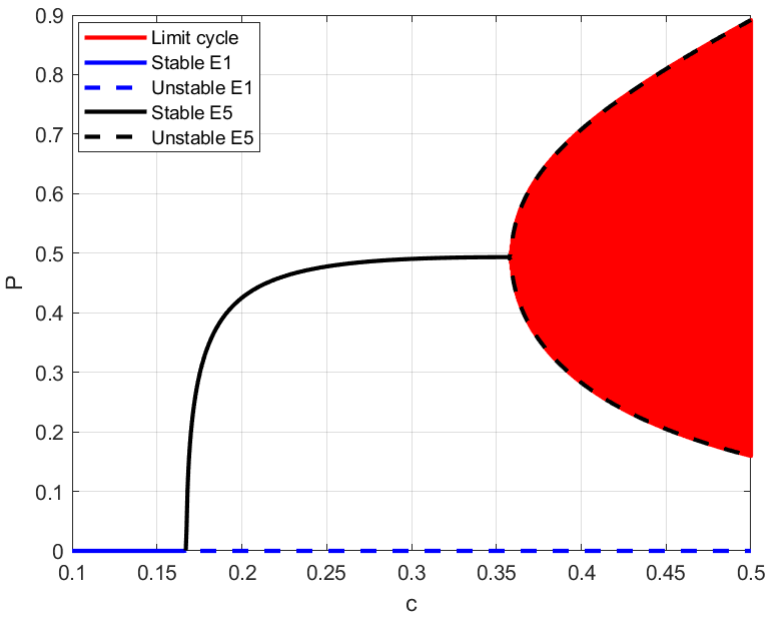}
			\caption{$ c-P $}
			\label{fig1b}
		\end{subfigure}
		\caption{\label{fig1} Bifurcation diagram of the system \eqref{GrindEQ__2_} with a weak Allee effect $ (h=0.2) $ and parameter values as in Table \ref{tab1}: (a) $N$ state and (b) $P$ state}
	\end{figure}

	The simulations in this section use the parameter values in Table \ref{tab1} and the predation conversion rate $c\in \left[0.1,\ 0.5\right]$. The impact of increasing the predation conversion rate on the convergence of the system's solution under the weak Allee effect condition is illustrated by the bifurcation diagram in Fig.\ref{fig1}.
	
	\begin{table}[h!]
		\caption{\label{tab1}Parameter values to observe the effect of the predation conversion rate}
		\begin{ruledtabular}
			\begin{tabular}{lccccccc}
				Parameter&$ r $&$ k $&$ w $&$ a $&$ b $&$ c $&$ \delta $\\
				\hline
				Value & 1 & 1 & 0.3 & 0.6 & 0.7 & $ 0.1 / 0.3 / 0.4 $ & 0.1
			\end{tabular}
		\end{ruledtabular}
	\end{table}
	
	\begin{figure}[h!]
		\centering
		\begin{subfigure}{0.49\textwidth}
			\centering
			\includegraphics[width=2.15in]{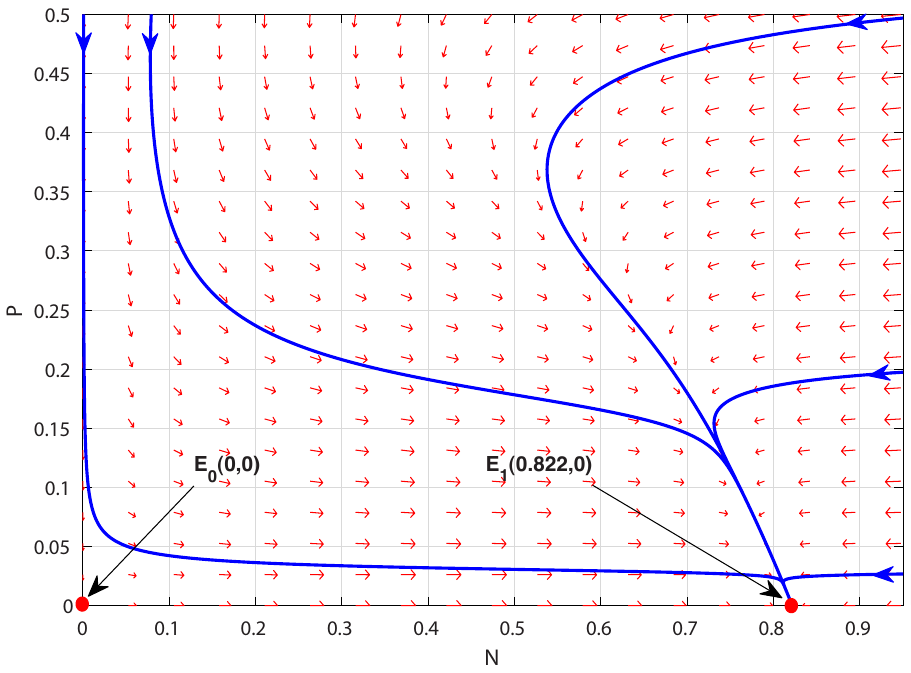}
			\caption{$ c=0.1 $}
			\label{fig2a}
		\end{subfigure}
		\hfill
		\begin{subfigure}{0.49\textwidth}
			\centering
			\includegraphics[width=2.15in]{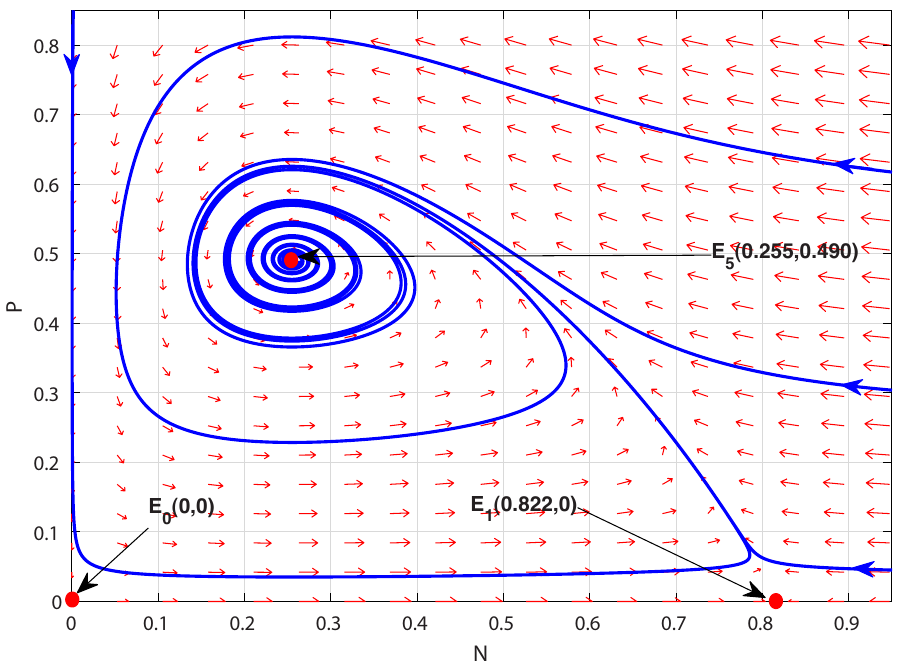}
			\caption{$ c=0.3 $}
			\label{fig2b}
		\end{subfigure}
		\\
		\begin{subfigure}{0.5\textwidth}
			\centering
			\includegraphics[width=2.15in]{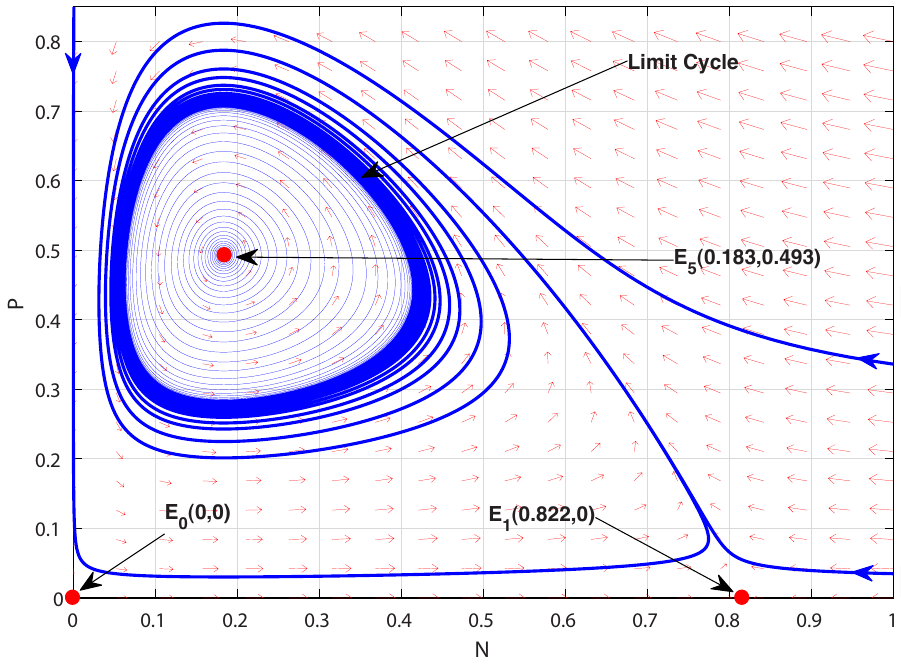}
			\caption{$ c=0.4 $}
			\label{fig2c}
		\end{subfigure}
		\caption{\label{fig2} Phase portraits of the system \eqref{GrindEQ__2_} with a weak Allee effect $ (h=0.2) $ and parameter values as in Table \ref{tab1}: (a)$ c=0.1$, (b) $c=0.3$, (c) $c=0.4$}
	\end{figure}
	
	The bifurcation diagram in Fig. \ref{fig1} identifies two bifurcation points with changes in the predation conversion rate: \( c_1^* = 0.167 \) and \( c_2^* = 0.359 \). A forward bifurcation shows a stability change of equilibrium point \(E_1 = (0.822,0)\) from locally asymptotically stable to unstable. For \(c < c_1^*\), \(E_1\) is locally asymptotically stable, but becomes unstable if \(c > c_1^*\). When \(E_1\) is unstable (\(c > c_1^*\)), the coexistence point \(E_5 = (0.452,0.425)\) is locally asymptotically stable for \(c_1^* < c < c_2^*\). If \(c > c_2^*\), a limit cycle emerges around \(E_5\), indicating a Hopf bifurcation. Fig. \ref{fig2} illustrates these stability changes in the phase portrait.
	
	\begin{figure}[h!]
		\centering
		\begin{subfigure}{0.49\textwidth}
			\centering
			\includegraphics[width=2.15in]{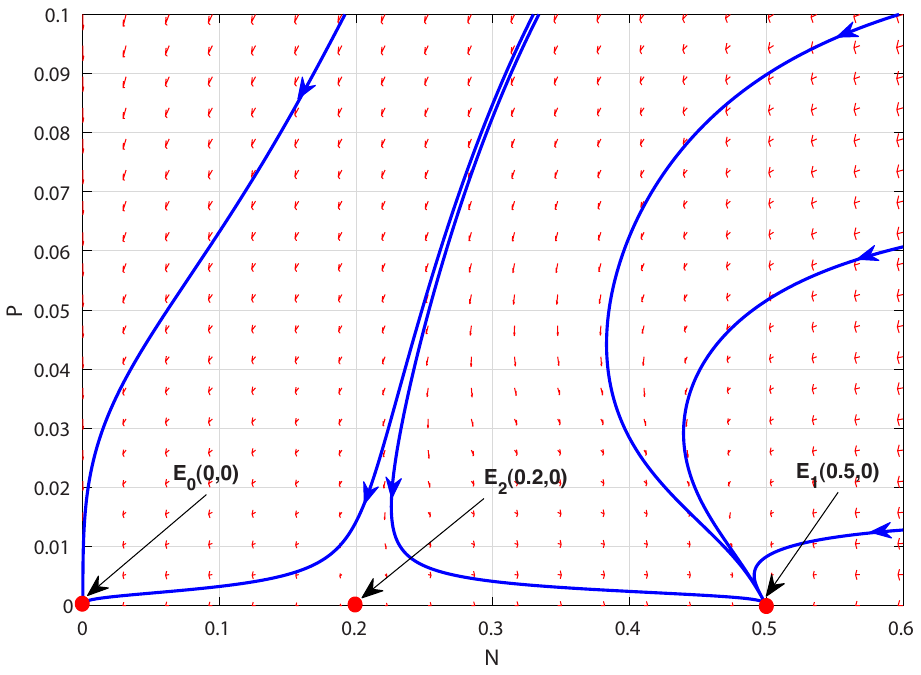}
			\caption{$ c=0.1 $}
			\label{fig4a}
		\end{subfigure}
		\hfil%
		\begin{subfigure}{0.49\textwidth}
			\centering
			\includegraphics[width=2.15in]{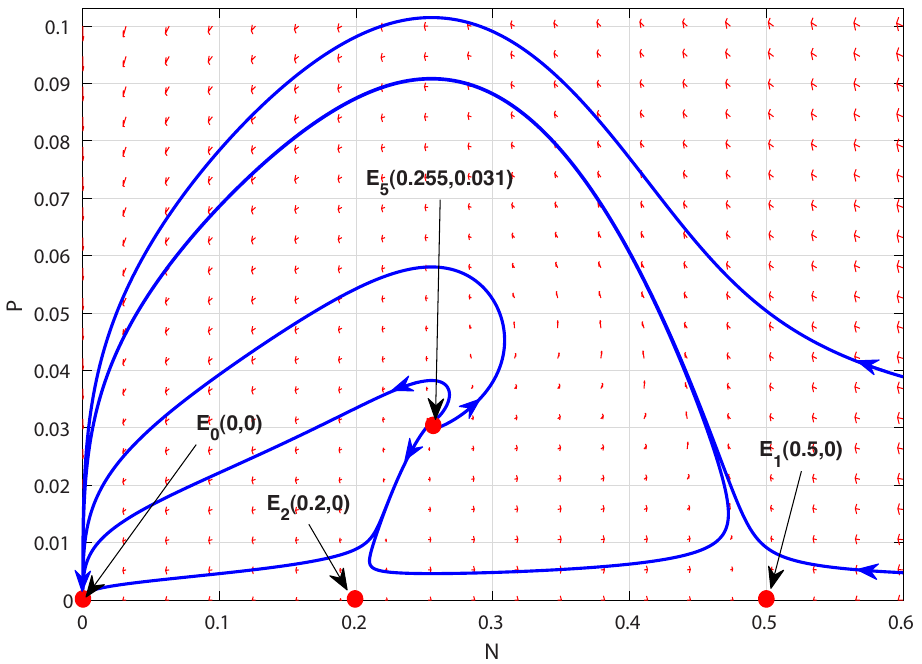}
			\caption{$ c=0.3 $}
			\label{fig4b}
		\end{subfigure}
		\\
		\begin{subfigure}{0.5\textwidth}
			\centering
			\includegraphics[width=2.15in]{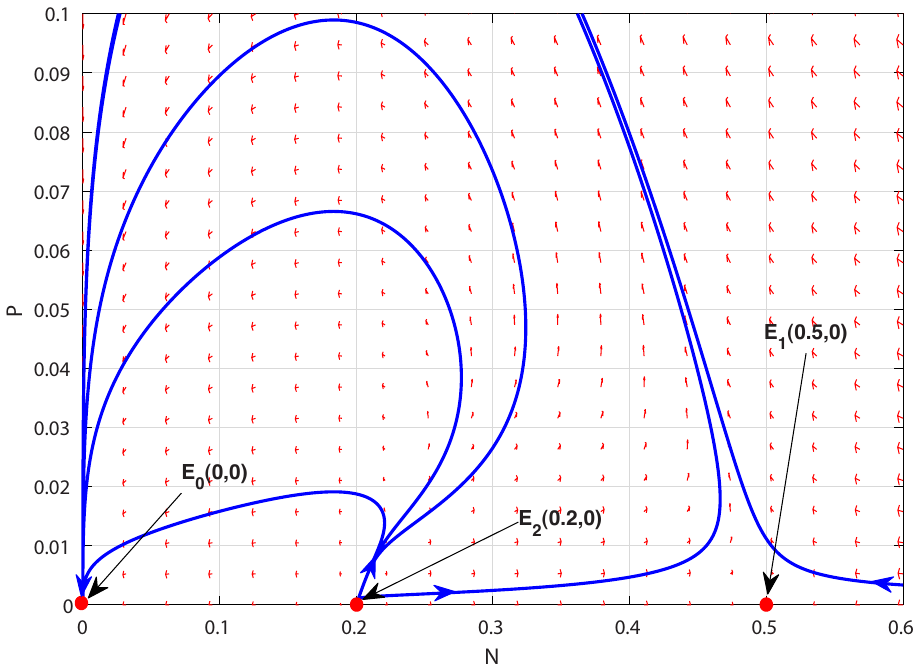}
			\caption{$ c=0.4 $}
			\label{fig4c}
		\end{subfigure}
		\caption{\label{fig4} Phase portraits of the system \eqref{GrindEQ__2_} with a strong Allee effect $ (h=0.4) $ and parameter values as in Table \ref{tab1}: (a)$\ c=0.1$, (b) $c=0.3$, (c) $c=0.4$}
	\end{figure}
	
	The phase portrait in Fig. \ref{fig2} illustrates the population dynamics of system \eqref{GrindEQ__2_} under a weak Allee effect ($h<w$). At a predation conversion rate of $c=0.1$ (Fig. \ref{fig2a}), the system has two equilibrium points: the trivial point $E_0=\left(0,0\right)$ and the axial point $E_1=\left(0.822,0\right)$, with convergence to $E_1$, indicating predator extinction and prey survival. Increasing the predation rate to $c=0.3$ (Fig. \ref{fig2b}) introduces a third equilibrium point, the coexistence point $E_5=(0.255,0.490)$, to which the system converges, suggesting stable coexistence of both populations. At $c=0.4$ (Fig. \ref{fig2c}), the system still has three equilibrium points but converges to a limit cycle around $E_5$, indicating long-term stable oscillations between prey and predator populations. Further simulations illustrate the effect of the predation conversion rate on population dynamics with a strong Allee effect in system \eqref{GrindEQ__2_}, with phase portraits given for varying $c$ values in Fig. \ref{fig4} and parameter values in Table \ref{tab1}.

	The phase portrait in Fig. \ref{fig4} illustrates the population dynamics in system \eqref{GrindEQ__2_} with a strong Allee effect ($h>w$). At a predation conversion rate of $c=0.1$ (Fig. \ref{fig4a}), the system has three equilibrium points: the trivial point $E_0=(0,0)$ and two axial points $E_1=(0.5,0)$ and $E_2=(0.2,0)$. The system converges to $E_0$ or $E_1$, indicating bistability, suggesting that prey survival depends on initial population size. At $c=0.3$ (Fig. \ref{fig4b}), there are four equilibrium points: the trivial point $E_0=(0,0)$, axial points $E_1=(0.5,0)$ and $E_2=(0.2,0)$, and the coexistence point $E_5=(0.255,0.031)$. The system converges to $E_0$, indicating extinction of both populations. At $c=0.4$ (Fig. \ref{fig4c}), the system has three equilibrium points: the trivial point $E_0=(0,0)$, and axial points $E_1=(0.5,0)$ and $E_2=(0.2,0)$. The coexistence point does not exist, and the system converges to $E_0$, suggesting that a higher predation conversion rate leads to the extinction of both populations.
	
	\subsection{The Impact of Environmental Protection Rate}
\begin{figure}[h!]
		%\centering
		\begin{subfigure}{0.49\textwidth}
			\centering
			\includegraphics[width=2.4in]{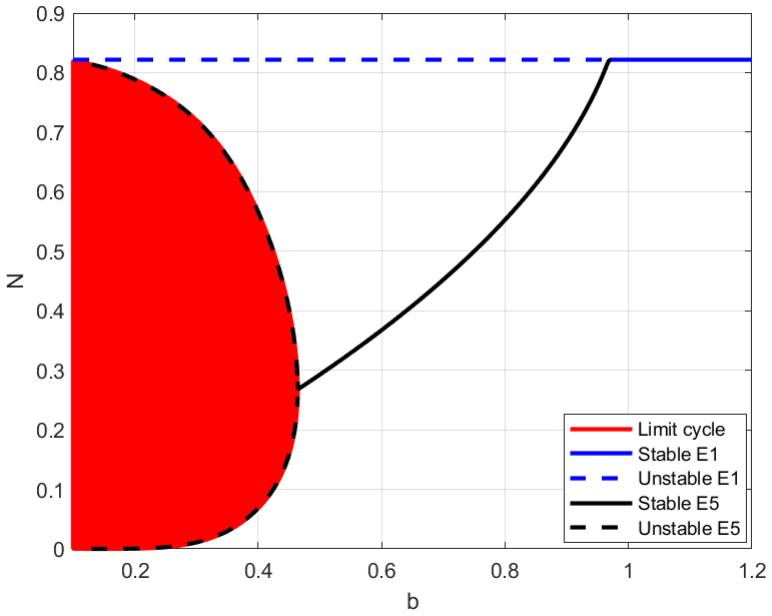}
			\caption{$ b-N $}
			\label{fig5a}
		\end{subfigure}
		\hfill
		\begin{subfigure}{0.49\textwidth}
			\centering
			\includegraphics[width=2.4in]{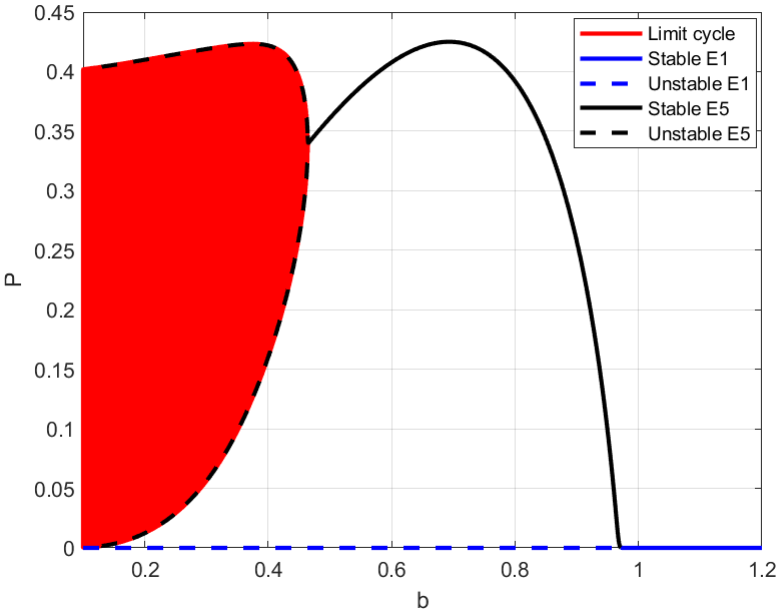}
			\caption{$ b-P $}
			\label{fig5b}
		\end{subfigure}
		\caption{\label{fig5} Bifurcation diagram of the system \eqref{GrindEQ__2_} with a weak Allee effect $ (h=2) $ and parameter values as in Table \ref{tab2}: (a) $N$ state and (b) $P$ state}
	\end{figure}

	The simulations in this section use the parameter values in Table \ref{tab2} and the environmental protection rate $b\in \left[0.1,\ 1.2\right]$. The impact of increasing the environmental protection rate on the convergence of the system's solution under the weak Allee effect condition is illustrated by the bifurcation diagram in Fig. \ref{fig5}.
	
	\begin{table}[h!]
		\caption{\label{tab2}Parameter values to observe the effect of the predation conversion rate}
		\begin{ruledtabular}
			\begin{tabular}{lccccccc}
				Parameter&$ r $&$ k $&$ w $&$ a $&$ b $&$ c $&$ \delta $\\
				\hline
				Value & 1 & 1 & 0.3 & 0.6 & $ 0.3 / 0.7 / 1.1 $ & 0.2 & 0.1
			\end{tabular}
		\end{ruledtabular}
	\end{table}
	
	\begin{figure}[h!]
		\centering
		\begin{subfigure}{0.49\textwidth}
			\centering
			\includegraphics[width=2.15in]{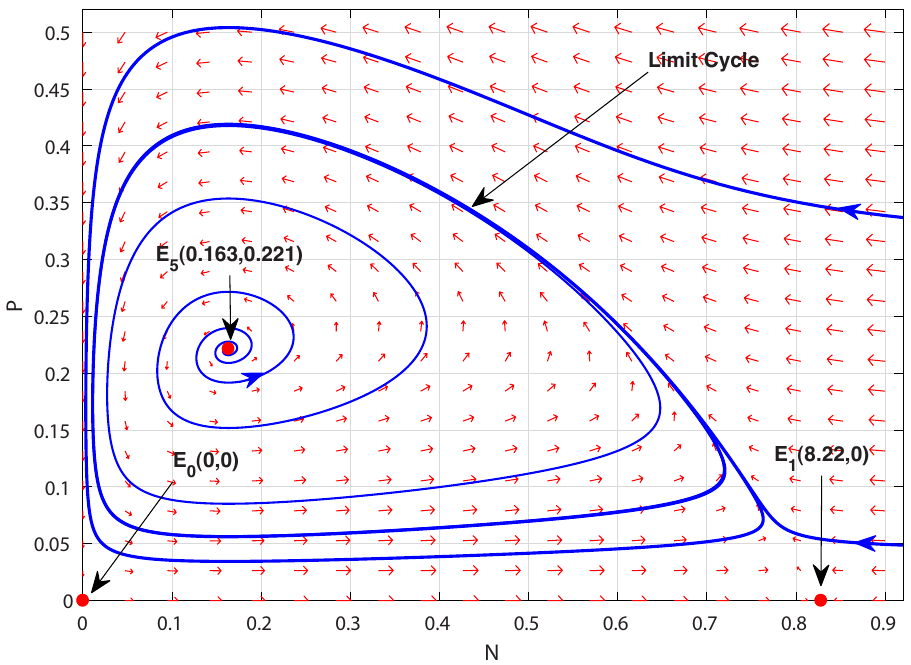}
			\caption{$ b=0.3 $}
			\label{fig6a}
		\end{subfigure}
		\hfill
		\begin{subfigure}{0.49\textwidth}
			\centering
			\includegraphics[width=2.15in]{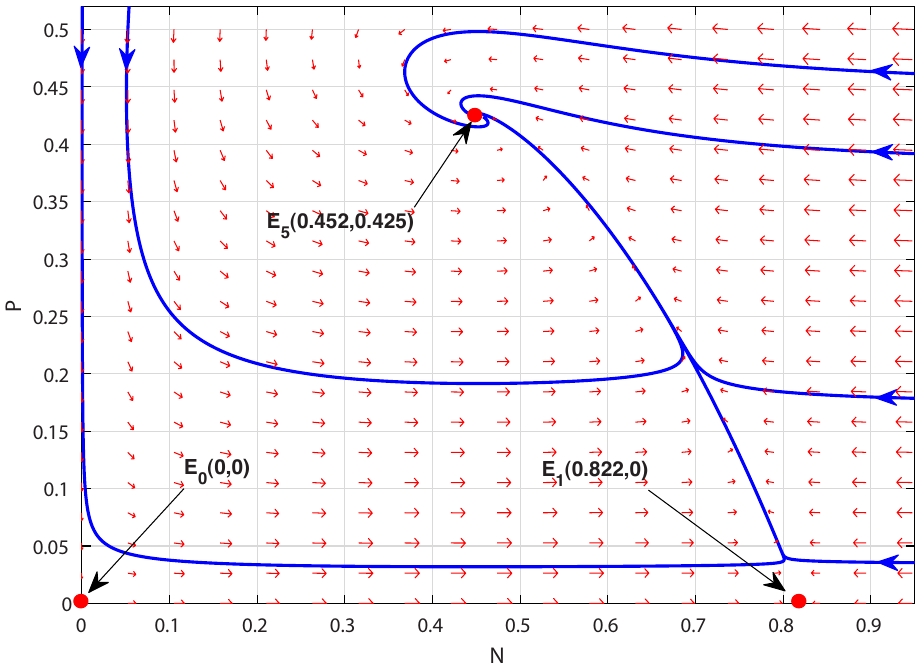}
			\caption{$ b=0.7 $}
			\label{fig6b}
		\end{subfigure}
		\\
		\begin{subfigure}{0.5\textwidth}
			\centering
			\includegraphics[width=2.15in]{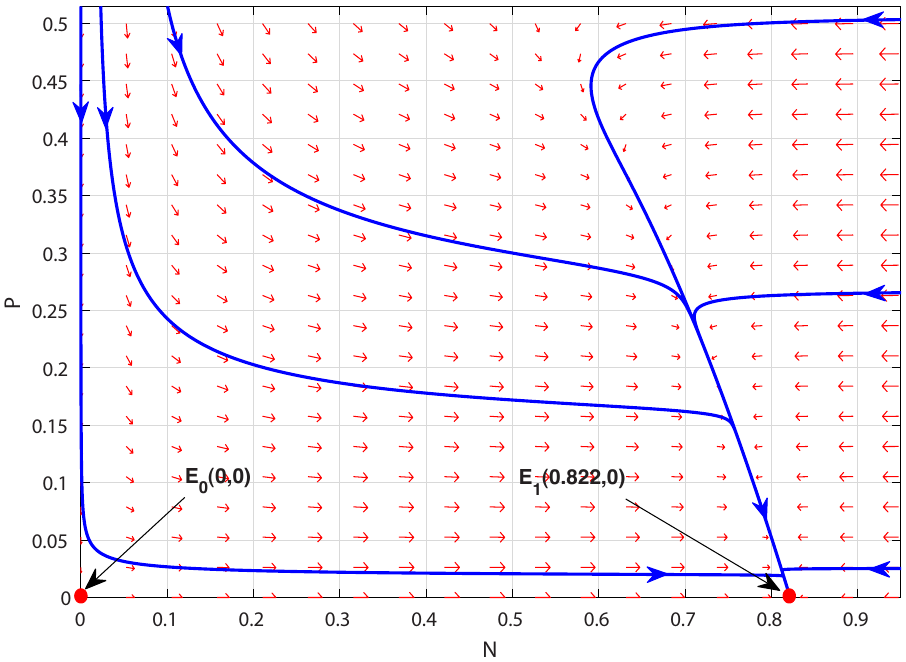}
			\caption{$ b=1.1 $}
			\label{fig6c}
		\end{subfigure}
		\caption{\label{fig6} Phase portraits of the system \eqref{GrindEQ__2_} with a weak Allee effect $ (h=0.2) $ and parameter values as in Table \ref{tab2}: (a)$\ b=0.3$, (b) $b=0.7$, (c) $b=1.1$}
	\end{figure}
	
	\begin{figure}[h!]
		\centering
		\begin{subfigure}{0.49\textwidth}
			\centering
			\includegraphics[width=2.15in]{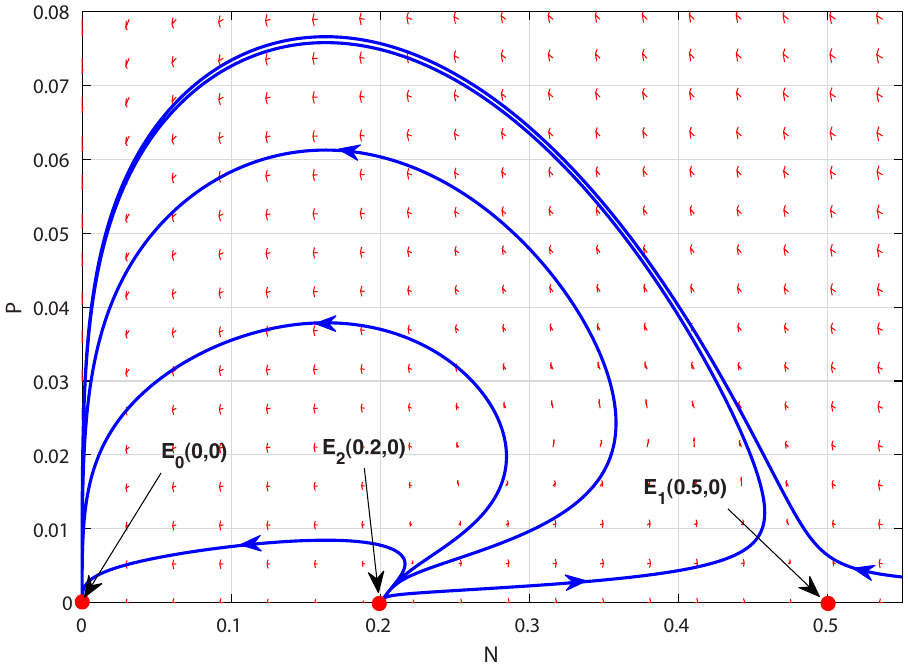}
			\caption{$ b=0.3 $}
			\label{fig7a}
		\end{subfigure}
		\hfill
		\begin{subfigure}{0.49\textwidth}
			\centering
			\includegraphics[width=2.15in]{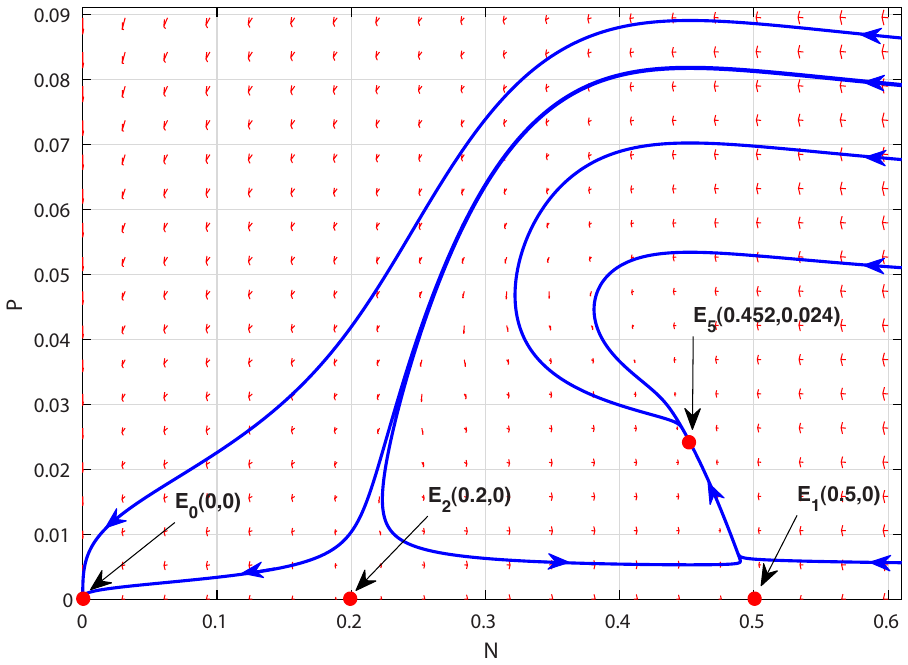}
			\caption{$ b=0.7 $}
			\label{fig7b}
		\end{subfigure}
		\\
		\begin{subfigure}{0.5\textwidth}
			\centering
			\includegraphics[width=2.15in]{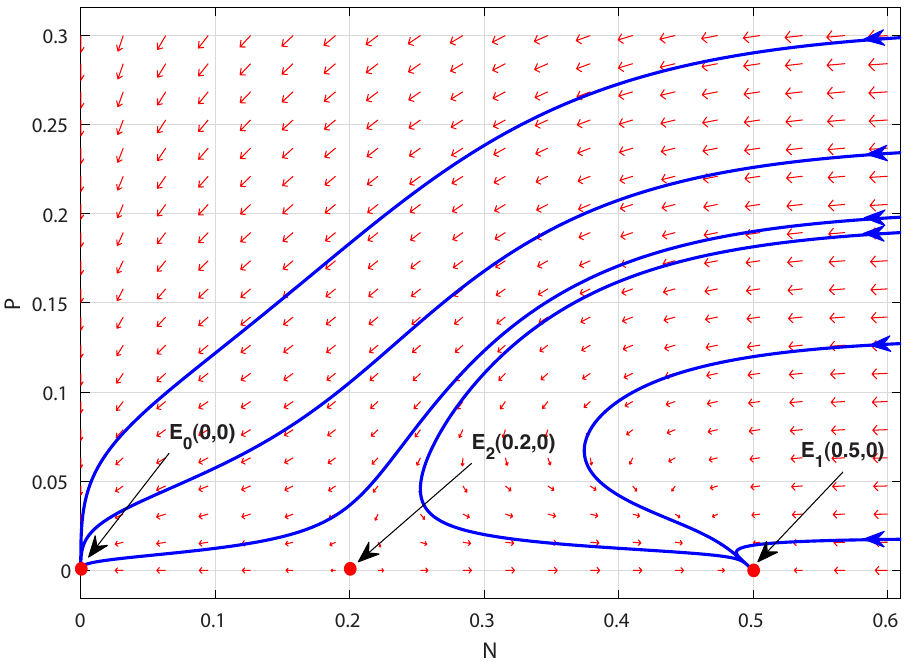}
			\caption{$ b=1.1 $}
			\label{fig7c}
		\end{subfigure}
		\caption{\label{fig7} Phase portraits of the system \eqref{GrindEQ__2_} with a strong Allee effect $ (h=0.4) $ and parameter values as in Table \ref{tab2}: (a)$\ b=0.3$, (b) $b=0.7$, (c) $b=1.1$}
	\end{figure}
	
	The bifurcation diagram in Fig. \ref{fig5} shows two bifurcation points for changes in the environmental protection rate: \( b_1^* = 0.465 \) and \( b_2^* = 0.972 \). The first bifurcation point marks a Hopf bifurcation, changing stability from a limit cycle to a stable coexistence point, \( E_5 \). For \( b < b_1^* \), a limit cycle occurs around \( E_5 \). When \( b_1^* < b < b_2^* \), \( E_5 \) is locally asymptotically stable, while \( E_1 = (N_1, 0) \) is unstable. If \( b > b_2^* \), a forward bifurcation changes the stability of \( E_5 \) to unstable and \( E_1 = (N_1, 0) \) to locally asymptotically stable. The phase portrait in Fig. \ref{fig6} illustrates these stability changes.

	The phase portrait in Fig. \ref{fig6} illustrates the population dynamics of system \eqref{GrindEQ__2_} under a weak Allee effect (\( h < w \)). At \( b = 0.3 \) (Fig. \ref{fig6a}), three equilibrium points exist: \( E_0 = (0,0) \), \( E_1 = (0.822,0) \), and \( E_5 = (0.163,0.221) \). The solution converges to a limit cycle around \( E_5 \), indicating stable oscillations. At \( b = 0.7 \) (Fig. \ref{fig6b}), the system still has three equilibrium points, with the solution converging to \( E_5 = (0.452,0.425) \), showing \( E_5 \) is locally asymptotically stable, while \( E_0 \) and \( E_1 \) are unstable, suggesting stable coexistence. At \( b = 1.1 \) (Fig. \ref{fig6c}), only \( E_0 \) and \( E_1 \) remain, with the solution converging to \( E_1 \), indicating predator extinction and prey survival. Additionally, simulations show the impact of the environmental protection rate on population dynamics with a strong Allee effect in system \eqref{GrindEQ__2_}. Phase portraits with parameter values from Table \ref{tab2} and variations in the environmental protection rate \( b \) are presented in Fig. \ref{fig7}.

	The phase portraits in Fig. \ref{fig7} illustrate population dynamics in system \eqref{GrindEQ__2_} under a strong Allee effect ($ h > w $). At an environmental protection rate of $ b = 0.3 $ (Fig. \ref{fig7a}), there are three equilibrium points: the trivial point $ E_0 = (0,0) $ and two axial points $ E_1 = (0.5,0) $ and $ E_2 = (0.2,0) $. Solutions converge to $ E_0 $, indicating potential extinction as $ E_0 $ is stable and the others are unstable. Increasing the protection rate to $ b = 0.7 $ (Fig. \ref{fig7b}) introduces a coexistence point $ E_5 = (0.452,0.024) $, with solutions converging to either $ E_0 $ or $ E_5 $, showing bistability where survival depends on initial population size. At $ b = 1.1 $ (Fig. \ref{fig7c}), only $ E_0 $ and $ E_1 $ remain, with solutions converging to these points, indicating potential extinction unless the initial prey population exceeds the Allee threshold.
	
	\section{Conclusion}
	This paper studies the dynamical analysis of predator-prey model incorporating the Allee effect and prey group defense. This model is proven to have non-negative solutions, existence, uniqueness, and boundedness, which confirm its validity in representing ecological phenomena. The model has three equilibrium points: the trivial point, the predator extinction point, and the coexistence point, all of which are locally asymptotically stable under certain conditions. If the Allee effect is weak, the trivial equilibrium point is unstable, while if the Allee effect is strong, the trivial equilibrium point is locally asymptotically stable. This implies that a strong Allee effect can lead to the extinction of both populations. Under weak Allee conditions, forward bifurcation and Hopf bifurcation occur at the predator extinction equilibrium point, while a strong Allee effect indicates bistability at both the trivial equilibrium point and the predator extinction equilibrium point. This means that prey can survive without the presence of predators, but a strong Allee effect can lead to prey extinction if the population size is very small. Numerical simulations supporting these findings are provided in the final section.
	
	%\nocite{*}
	\bibliography{aipsamp}% Produces the bibliography via BibTeX.
	
\end{document}